\documentclass[12pt]{article}
\usepackage{amssymb, amsmath, amsthm, amsfonts}
\usepackage[italian, english]{babel}
\usepackage{graphicx}
\usepackage[T1]{fontenc}
\usepackage[latin1]{inputenc}
\usepackage{times}

\usepackage{epsfig}
\usepackage{graphicx}
\usepackage{hyperref}

\newtheorem{theorem}{Theorem}[section]
\newtheorem*{theorem*}{Theorem}
\newtheorem{lemma}[theorem]{Lemma}
\newtheorem{definition}[theorem]{Definition}
\newtheorem{corollary}[theorem]{Corollary}
\newtheorem{proposition}[theorem]{Proposition}
\newtheorem{remark}[theorem]{Remark}

\newcommand {\E}{\mathbb E}
\newcommand{\R}{\mathbb R}
\newcommand{\N}{{\mathbb N}}
\newcommand{\sfe}{{\mathbb S}^1}
\newcommand{\cl}{{\rm cl}}
\newcommand{\Tan}{{\rm Tan}}
\newcommand{\interno}{{\rm int}}
\newcommand{\co}{{\rm co}}
\newcommand{\nor}{{\rm Nor}}
\newcommand{\reach}{{\rm reach}}
\newcommand{\dist}{{\rm dist}}
\newcommand{\diam}{{\rm diam}}
\newcommand{\length}{{\rm length}}
\newcommand{\arc}{{\rm arc}}
\newcommand{\meas}{{\rm meas}}
\newcommand{\per}{{\rm per}}

\begin{document}
	
\selectlanguage{english}
\title{Plane $R$-paths and their rectifiability property}
\author{Nico Lombardi   \and  Marco Longinetti \and Paolo Manselli \and Adriana Venturi}
\date{}
\maketitle

\begin{abstract}
	
	A family of plane oriented continuous paths depending on a fixed real positive number $R$ is considered. For any point $x$ on the path, the previous points lie out of any circle of radius $R$ having at $x$ interior normal in a suitable tangent cone to the path at $x$. These paths are locally descent curves of a nested family sets of reach $R$.
	
	Avoiding any smoothness requirements, we get angle estimate and not intersection property. Afterwards we are able to estimate the lenght and detour of this curve.
	
	\medskip
	{\noindent Keywords: Steepest descent curves, sets with positive reach, length of curves, detour.}
	
	{\noindent 2010 AMS subject classification: Primary: 28A75, 52A30;  Secondary: 34A26, 34A60.}
	
\end{abstract}

\section{\textbf{Introduction}}\label{intro}

\noindent
Let $g \colon [a,b] \to \R^2$ be a continuous mapping, oriented according to the increasing variable. For $t\in [a,b]$, $x=g(t)$, let
$$g_x:=g([a,t]),\quad g^x:=g([t,b]).$$

\noindent
Let us assume that for every $x\in g([a,b])$ there is a constraint $Q_x$ on $g$; then what extra properties does the mapping $g$ satisfy?

Different constraints for the previous points have been considered by several authors and for each constraint different classes of family of curves have been studied. Most of these curves are related to descent curves, as in \cite{Daniilidis/Deville} and \cite{Mainik}.

\noindent
Let us list some interesting examples:

\noindent
\textbf{Example 1}. Let the constraint $Q_x$ be that, for $a\leq t_1 <t_2<t$ and $x=g(t)$, 
\begin{equation}\label{first inequality}
|g(t_1)-g(t_2)|\leq |g(t_1)-g(t)|.
\end{equation}
\noindent
This family of mappings has been studied in \cite{Manselli-Pucci} as steepest descent curves of quasi-convex functions and
they have been studied in several dimensions and called self expanding curves in \cite{MLV}.

\noindent
In \cite{MLV} it was proved that $g$ is rectifiable: as a consequence, if $s\in [0,L]\longmapsto x(s)$ is the representation of $g$ in arc length, differentiating \eqref{first inequality}, it turns out that

\begin{equation}\label{g intro}
g^{x(s)} \subseteq \{ y\in \R^2:\ \langle y-x(s),x'(s) \rangle \leq 0  \}.
\end{equation}
This means that the previous points to $x(s)$ lie out of the half-space orthogonal to $x'(s)$ at $x$. 

\noindent
It turns out that the continuous mapping $g$ is an injective rectifiable curve and there exists $c>0$ so that:
\begin{equation}\label{Length intro}
L=\length(g)\leq c\ |g(b)-g(a)|.
\end{equation}
Moreover there exists a nested family of convex sets $\Omega_s$ such that $x(s)\in \partial \Omega_s$ and $x'(s)$ belongs to the normal cone to $\Omega_s$
at $x(s)$, see \cite{MLV}. If $f$ is a regular quasi-convex function with level sets $\partial \Omega_s$, then $s \longmapsto x(s)$ is a steepest descent curve for $f$.

Let us remark that continuous mappings $g$ defined as in inequality \eqref{first inequality}, but with opposite order, were called self-approaching curves in 
\cite{Langetepe}, self-contracting curves in \cite{Daniilidis} and curves of descent in \cite{Mainik}.

\vspace{0.25cm}
\noindent
\textbf{Example 2}. Let $\phi\in (0,\pi)$ and let us assume, as constraint $Q_x$, the existence of a wedge
of opening $\phi$ and vertex $x$ containing $g^x$.
These mappings $g$ are called $\phi$-self-approaching curves in  \cite{Langetepe2}, 
see  Definition \ref{defphiself}.
Then $g$ is an injective rectifiable curve and \eqref{Length intro} holds for a best constant $c$ esplicitly computed.

\vspace{0.25cm}
\noindent
\textbf{Example 3}. Let $\lambda\in [-1,1)$ and $\alpha=\arccos(\lambda)$.

Let $\gamma$ be a plane continuous path with nonzero right derivatives $\gamma'(t)$ at $x=\gamma(t)$. Let the constraint $Q_x$ be:
\begin{equation*}
<\gamma'(t),\gamma(u)-\gamma(t)> \leq \lambda |\gamma'(t)| |\gamma(u)-\gamma(t)|
\end{equation*}
for every $u < t$; geometrically the constraint is that $\gamma_x$ does not intersect the open cone with axis $\gamma'(t)$ and opening $2\alpha$; these curves have been called $\lambda$-eels in \cite{Daniilidis/Deville}.

In \cite[Theorem 6.1]{Daniilidis/Deville} it is proved that bounded $\lambda$-eels have finite length.

\vspace{0.25cm}
\noindent
\textbf{Example 4}. Let us assume that there exists a family of nested increasing convex sets $\Gamma_t$, so that $x=g(t)\in \Gamma_t$, $0\leq t\leq 1$; for every $x=g(t)$, $0<t<1$, let the constraint $Q_x$ be the following: for every unit vector $v$ in a suitable tangent set to $g^x$ at $x$, there exists a support hyperplane $P(x)$ to the set $\Gamma_t$, at $x=g(t)$, such that the angle between $P(x)$ and $v$ is at least $\frac{\pi}{2}-\theta$. Then in \cite{Mainik2} it is proved that \eqref{Length intro} holds for $\theta$ sufficiently small.

\vspace{0.25cm}
\noindent
\textbf{Example 5}. Let us assume that $R>0$ is fixed and $\gamma$ a rectifiable curve.
Let $[0,L]\ni s \longmapsto  x(s)\in  \R^2$ be its representation in arc length.
Let the constraint $Q_x$ be such that $\gamma([0,s])$ has empty intersection with
the open disk of radius $R$ and centered in $x(s)+Rx'(s)$, denoted by $B(x(s)+Rx'(s))$.
In \cite{MLV2}  and \cite{steepest I}, $\gamma$ has been called an \textit{$R$-curve}, see also Definition \ref{defRcurve}.

\noindent
In \cite{steepest I} and \cite{MLV2} it is proved that, if $\gamma([0,L])$ is contained in a small ball, 
then \eqref{Length intro} holds. Moreover if $\gamma([0,L])$ is contained in a bounded set $G$, 
then \eqref{Length intro} holds in a weaker form: there exists a constant $c(\diam (\gamma))>0$ such that
\begin{equation*}
\length(\gamma)\leq c( \diam(\gamma))\ |\gamma(0)-\gamma(L)|.
\end{equation*}

\vspace{0.25cm}
Let us notice that in all the Examples $1$, $2$, $3$ and $4$, the rectifiability hypothesis was not needed.
It is natural to ask if, in Example $5$, the rectifiability hypothesis could be also dropped. That is the main goal of the present work. Moreover it is not required any regularity assumption on the curve.

Let us consider the path satisfying the following constraint:
$g\in [a,b] \to \R^2$ is a continuous mapping and for every $t\in (a,b]$,
$$g([a,t])\cap B(g(t)+Rl)= \varnothing,$$
where $l$ is a unit vector in a suitable tangent set to $g$ at $g(t)$, see Definition \ref{def R-path}; $g$ will be called \textit{$R$-path}. The name \textit{path} is used for a continuous, not necessarily rectifiable map; the name \textit{curve} for rectifiable map.

\noindent
In Theorem \ref{First Main Thm} it will be proved that every injective $R$-path in a small ball is rectifiable.

\noindent
In Theorem \ref{Second Main Thm} it will be proved that every $R$-path in a small ball is injective.

Then every $R$-path turns out to be a rectifiable $R$-curve (Theorem \ref{R-path and curve}) and, in an open ball of radius $R$, the $R$-paths are descent curves of a nested family of sets of reach $\geq R$, Theorem \ref{R-path in a disk}.

\section{\textbf{Definitions and preliminaries}}\label{preliminaires}

\noindent
$R$ will be a fixed positive real number throught this work. Let 
$$ B(z,\rho)=\{ x\in \R^2\,:\,|x-z|<\rho \},\ \quad\
\sfe=\partial B(0,1)$$
and let $D(z,\rho)$  be the closure of $B(z,\rho)$.
The notations $B_\rho, D_\rho$ will also be used for  balls of radius $\rho$, if no ambiguity arises for their center.
If a ball is written with a center only, then the radius  will be $R$.
The usual scalar product between vectors $u,v\in \R^2$ will be denoted by $\langle u,v\rangle $ and by $|u|$ the Euclidean norm in $\R^2$.

Let $K \subset \R^2$;
$\interno(K)$ will be  the interior of $K$,  $\partial K$ 
the boundary of $K$, $\cl(K)$ the closure of $K$,
$K^c=\R^2 \setminus K$ and $\per(K)$ will be the perimeter of $K$, if $K$ is convex. For every set $S \subset \R^2$, $\co(S)$ is the convex hull of $S$.

Let $K$ be a non empty closed set.
Let $q\in  K$; the {\em tangent cone} of $K$ at $q$  is defined as
\begin{equation*}\label{tangent}
\Tan(K,q)=\{v \in \R^2: \forall \varepsilon > 0 \; \exists \, x\in K,\; r > 0\ \textnormal{with}\ |x-q|< \varepsilon,\; |r(x-q)-v|< \varepsilon\}.
\end{equation*}

\noindent 
Let us recall that, if $q$ is not an isolated point of $K$, then
$$\sfe\cap \Tan(K,q)= \bigcap_{\varepsilon >0} \cl(\{\frac{x-q}{|x-q|},\, q\neq x\in K\cap B(q,\varepsilon)\}).$$


\noindent
The {\em normal cone} at
$q$ to $K$ is the non empty closed convex cone,  given by:
\begin{equation}\label{normalcone} \nor(K,q)=\{u\in\R^2: \langle u,v\rangle \le 0, \quad
\forall v \in \Tan(K,q)\}.
\end{equation}
When $q \in \interno(K)$, then $\Tan(K,q)=\R^2$ and  $\nor(K,q)$ reduces to zero.
In two dimensions cones will be called angles  with vertex 0.



\noindent
Let $b \in \R^2\setminus A$; $b$ has a unique projection point onto $A$ if there exists a unique point $a\in A$ satisfying
$|b-a|=\dist(b,A)$.

Let $A$ be a closed set. If $a\in A$ then $\reach(A,a)$ is the supremum of all numbers $\rho$ such that every $x\in B(a,\rho)$ has a unique projection point 
onto $A$ and
$$\reach(A):=\inf \{\reach(A,a): a \in A \},$$
see \cite{Federer}.






\noindent
Let $x,y \in \R^2, |x-y| < 2R.$ Let us define
\begin{equation}\label{deffusotto}
\mathfrak{h}(x,y,R)=\bigcap \{D_R:\ D_R\supset \{x,y\}\}.
\end{equation}

\begin{proposition}\cite{Ratay1},\cite[Theorem 3.8]{Colman} \label{fusotto}
	A closed set $A$ has reach equal or greater than $R$ if and  only if $\forall x,y \in A$ with $|x-y|< 2R$ the set 
	$A\cap  \mathfrak{h}(x,y,R)$ is connected. 
\end{proposition}

\begin{definition}\label{R-hull}
	Given $A$ a closed set in $\R^2$, let us define $\co_R(A)$, the {\em  $R$-hull } of $A$, as the closed set containing $A$, such that
	\renewcommand{\theenumi}{\roman{enumi}}
	\begin{enumerate}
		\item $\co_R(A)$ has reach greater or equal than $R$;
		\item if a set $B \supseteq A $ and $\reach(B)\geq R$, then $B\supseteq \co_R(A)$.
	\end{enumerate}
\end{definition}
See \cite[pp.105-107]{Colman}   for the properties of $R$-hull. The $R$-hull of a closed set $A$ may not exist, see \cite[Remark 4.9]{Colman}. However if $\co_R(A)$ exists, it can be shown that
$$   \co_R(A)=\bigcap\{ B_R^c:\ B_R\cap A=\emptyset \}. $$

\noindent
Moreover a sufficient condition for the existence of the $R$-hull has been proved:

\begin{proposition}\cite[Theorem 4.8]{Colman}\label{propRhull} 
	If $A$ is a closed connected subset of an open circle of radius $R$ in $\R^2$, then $A$ has $R$-hull.
\end{proposition}

\begin{remark}\label{R-hull inclusion}
	\textnormal{Let $A\subset B$ be two closed connected subsets of an open circle of radius $R$, then } $$\co_R(A)\subset \co_R(B).$$
\end{remark}

\begin{definition}\label{defRcurve}
	An $R$-curve $\gamma \subset \R^2$ is a rectifiable oriented curve with arc length parameter
	
	\noindent $s\in [0,L]$, tangent vector $\textbf{t}(s)=x'(s)$ such that the inequality
	\begin{equation}\label{Rcurve rec}
	|x(s_1)-x(s)-R\ \textbf{t}(s)|\geq R
	\end{equation} 
	holds for almost all $s$ and for $0\leq s_1\leq s \leq L$.
\end{definition}

Let us notice the following equivalent formulations of \eqref{Rcurve rec} for $0\leq s_1<s\leq L$:
\begin{align}
&|x(s_1)-x(s)|^2 -2R\langle x(s_1)-x(s), \textbf{t}(s) \rangle \geq 0;\label{dis with tagent vector 1}
\\ & \langle x(s)-x(s_1), \textbf{t}(s) \rangle \geq -\frac{|x(s_1)-x(s_2)|^2}{2R};\label{dis with tagent vector 2}
\\ & \langle \frac{x(s_1)-x(s)}{|x(s_1)-x(s)|}, \textbf{t}(s) \rangle \leq \frac{|x(s_1)-x(s)|}{2R},
\quad \textnormal{if}\quad x(s_1)\neq x(s).\label{dis frazione}
\end{align}

The next property of $R$-curves is proved in \cite{MLV}.

\begin{proposition}\cite[Corollary 3.3]{MLV}\label{property R-curve}
	An $R$-curve does not intersect itself.
\end{proposition}

\section{\textbf{$R$-paths and their properties}}

\noindent
Let $g\colon [a,b] \to \R^2$ be a continuous mapping satisfying $g(a)\neq g(b)$.

Let $x=g(t)$, $x_1=g(t_1)$ and $x_2=g(t_2)$, with $a\leq t \leq b$ and $a\leq t_1<t_2\leq b$; let: $$g_x:=g([a,t]),\quad g_{x_1x_2}:=g([t_1,t_2]),\quad g^x:=g([t,b]).$$

If $A\subset \R^2$,
\begin{equation*}
\stackrel{-1}{g}(A)=\{ t\in [a,b]:\ g(t)\in A \}
\end{equation*}
will be the preimage of $A$ under $g$.

\noindent
Let $t\in [a,b]$, $x=g(t)$. Let us define:
$$\Theta_{-}(x)=-\Tan(g_x,x),\quad \Theta_{+}(x)=\Tan(g^x,x).$$
\noindent
If $x$ is not a starting point, then $\Theta_{-}(x)\supsetneqq \{ 0 \}$. If $x$ is not an end point, then $\Theta_{+}(x)\supsetneqq \{ 0 \}$.

Let us remark that $\Theta_{-}(x)$ and $\Theta_{+}(x)$ exist for any $x\in g([a,b])$ without any regularity assumption, since they are the two sets of all left-limits of chords and all right-limits of chords.

\begin{definition}\label{def R-path}
	The mapping $g\colon [a,b] \to \R^2$ will be called an $R$-path if it 
	is a continuous mapping satisfying $g(a)\neq g(b)$ and, for all $x=g(t)$, $a\leq t \leq  b$, the following condition: for every $l\in [\Theta_{-}(x)\cup\Theta_{+}(x)]\cap \sfe$, $$g_x \subset B^c(x+Rl)$$ holds.
\end{definition}

\begin{remark}\label{equivalent R-path}
	\textnormal{The previous inclusion can be written as:}
	\begin{equation}\label{Condition R-path}
	|g(\tau)-g(t)|^2-2R\langle g(\tau)-g(t), l \rangle \geq 0,\quad a\leq \tau\leq t\leq b,
	\end{equation}
	\textnormal{which extends inequality \eqref{dis with tagent vector 1} and the equivalent inequalities \eqref{dis with tagent vector 2} and \eqref{dis frazione} to an $R$-path. Furthermore let us observe that a rectifiable $R$-path is an $R$-curve.}
\end{remark}

\begin{remark}\label{Joinying curve}
	\textnormal{Let $g\colon [a,b] \to \R^2$ be an $R$-path and let $t(\tau)\colon [\alpha,\beta]\to [a,b]$ be a strictly monotonic continuous function with inverse function $\tau(t)$. Then the function
		$$h\colon [\alpha,\beta]\to \R^2,$$ defined as $$h(\tau)=g(t(\tau)),\quad \forall\ \tau\in[\alpha,\beta],$$ is an $R$-path with $h([\alpha, \beta])=g([a,b])$.}
	
	\noindent
	\textnormal{Let us notice that we can reparametrize an $R$-path in such way that $\stackrel{-1}{g}(\{x\})$ has no interior points for every $x$.}
\end{remark}

\begin{lemma}\label{joining curve}
	Let $\gamma_i\colon [0,L_i]\to \R^2$ be $R$-curves in arc length representation $x=x_i(s)$, $0\leq s\leq L_i$ and $i=1,2$. Assume that 
	\begin{itemize}
		\item $x_1(L_1)=x_2(0)$,
		\item for almost all $s_2 \in [0,L_2]$
		\begin{equation}\label{inclusion joint lemma}
		\gamma_1([0,L_1])\subseteq B^c(x_2(s_2)+Rx'_2(s_2)).
		\end{equation}
	\end{itemize}
	Let $\gamma\colon [0,L_1+L_2]\to \R^2$ be defined as
	\begin{equation*}
	\gamma(s)=\begin{cases}
	x_1(s), & 0\leq s\leq L_1,\\
	x_2(s-L_1),& L_1< s\leq L_1+L_2.
	\end{cases}
	\end{equation*}
	Then $\gamma$ is an injective rectifiable $R$-curve and
	$$\gamma_1([0,L_1])\cup \gamma_2([0,L_2])=\gamma([0,L_1+L_2]).$$
\end{lemma}

\noindent
{\it Proof.} The function $\gamma\colon [0,L_1+L_2]\to \R^2$ is Lipschitz continuous; for $s,s_1 \in [0,L_1]$, such that $s_1<s <L_1$, \eqref{Rcurve rec} holds for $\gamma=\gamma_1$.
Similarly if $L_1 < s_1 <s< L_1+L_2$.

If $s_1\in [0,L_1]$ and $s_2\in [L_1,L_1+L_2]$, then \eqref{Rcurve rec} follows from \eqref{inclusion joint lemma}. So $\gamma$ is a rectifiable $R$-curve and it is injective, by Proposition \ref{property R-curve}.
\qed

\begin{lemma}
	Let $g\colon [a,b]\to \R^2$ be an $R$-path, $t_0\in (a,b)$ and $x=g(t_0)$. Let $u\in \sfe$ satisfying $B(x+Ru)\cap g([a,b])=\varnothing$. Then $$<l,u>\ \geq 0\quad \forall\ l\in \Theta_{-}(x)$$ and $$<l,u>\ \leq 0\quad \forall\ l\in \Theta_{+}(x)$$ hold.
\end{lemma}

\noindent
{\it Proof.} Let $l\in \Theta_{-}(x)$ (similar proof for $l\in \Theta_{+}(x)$). Then $-l\in \Tan(g_x,x)$ and there exists $\{ t_n \}_{n\in \N}\subseteq (a,t_0)$ so that $x\neq g(t_n)\rightarrow x$ and $$\frac{g(t_n)-x}{|g(t_n)-x|}\rightarrow -l.$$

\noindent
As $g(t_n)\in B^c(x+Ru)$, $|g(t_n)-(x+Ru)|^2\geq R^2$ implies
$$<\frac{g(t_n)-x}{|g(t_n)-x|},u>\ \leq \frac{|g(t_n)-x|}{2R}.$$
As $n\rightarrow +\infty$, the inequality $<-l,u>\leq 0$ follows.
\qed

\begin{lemma}\label{Successione}
	Let $g\colon [a,b]\to \R^2$ be an $R$-path, $x=g(t)\in g([a,b])$, $x\neq g(a)$.
	
	Let us assume that
	\begin{itemize}
		\item there exists a sequence $\{ u_k \}_{k\in \N}\subseteq \sfe$ and $u\in \sfe$ such that $u_k \rightarrow u$;
		\item there exists a sequence $\{ t^{(k)} \}_{k\in \N}\subseteq [a,b]$ satisfying $ g(t^{(k)}) \rightarrow x$, with $x\neq g(t^{(k)})$, $t^{(k)}<t$, for every $k\in \N$, and $g(t^{(k)})\in \partial B(x+Ru_k)$.
	\end{itemize}
	Then there exists $\overline{u}\in \Tan(g_x,x)$ satisfying $<\overline{u},u>\ =0$.  
\end{lemma}

\noindent
{\it Proof.} As $g(t^{(k)})\in \partial B(x+Ru_k)$, then 
\begin{equation}\label{Equality Lemma}
<\frac{g(t^{(k)})-x}{|g(t^{(k)})-x|},u_k>=\frac{|g(t^{(k)})-x|}{2R}
\end{equation}
holds.

\noindent
Let, up to subsequences, $\overline{u}:=\lim_{k \rightarrow +\infty} \frac{g(t^{(k)})-x}{|g(t^{(k)})-x|}$. 

As $k\rightarrow +\infty$, $<\overline{u},u>=0$ follows from \eqref{Equality Lemma}.
\qed

\begin{theorem}\label{4.1lavoro1}Let $g\colon [a,b]\to \R^2$ be an $R$-path such that $g([a,b]) \subset B_R$. Let $x\in g([a,b])$ and let
	\begin{equation}\label{defwx}
	W_x=\{ u\in \sfe : g_x \subseteq B^c(x+Ru) \}.
	\end{equation}
	Then
	\begin{equation}\label{U+UnioneU-inclusoWx}
	(\Theta_{-}(x)\cup\Theta_{+}(x))\cap \sfe \subset W_x.
	\end{equation}
	Moreover
	\begin{equation}\label{wx=N}
	W_x=\nor(\co_R(g_x),x)\cap \sfe.
	\end{equation}
\end{theorem}

\noindent
{\it Proof.} \eqref{U+UnioneU-inclusoWx} follows immediately from Definition \ref{def R-path}. The proof of \eqref{wx=N} is the same as in \cite[Lemma 4.3]{steepest I} by changing $\gamma_x$ into $g_x$.

Let us notice that the assumption of $g([a,b]) \subset B_R$ implies the existence of $\co_R(g_x)$, by Proposition \eqref{propRhull}.
\qed

\section{\textbf{Injective $R$-paths}}
Let us first recall some geometric definitions and properties, see \cite{steepest I}.
\begin{definition}\label{defangcur} Let $b$ and $c$ two  distinct points in the plane with $|b-c|< 2R$.
	Let us consider the open circles $B_R(b)$ and $B_R(c)$ and $x\in \partial B_R(b) \cap \partial B_R(c)$. Let $l$ be the line through $b$ and $c$ and $ H $ be the half plane with boundary $l$ containing $x$. 
	
	The unbounded region $\check{ang}(bxc):\equiv
	B_R^c(b) \cap B_R^c(c)\cap H$
	will be   called   {\em curved   angle}.
	Moreover 
	$$\meas (\check{ang}(bxc)):=\meas(\Tan(\check{ang}(bxc),x)\cap \sfe)$$
	is the measure of the angle between the half tangent lines at $x$ to the boundary of $\check{ang}(bxc)$. 
\end{definition}

When $x,y$ are points on a circumference $\partial B$, let us denote by $\arc(x,y)$ the shorter arc on $\partial B$  from $x$ to $y$.

\begin{proposition} \cite[Lemma 4.4]{steepest I}\label{geometrico} 
	Let $x, x_2\in \R^2$, $|x-x_2|< R$. 
	Let $B^2:=B_R(b)$, with $b\in \R^2$ and $x,x_2$ on $\partial B^2$. 
	Let $B^*:=B_R(c_*)$ the ball orthogonal at $x_2$ to $ B^2$ such that $x\in B^*$. Let us assume that there exists $x_1\in (B^*\cup B^2)^c$ with the properties (see Fig.1 in \cite{steepest I}):
	\renewcommand{\theenumi}{\roman{enumi}}
	\begin{enumerate}
		\item $|x_1-x| < R, |x_2-x_1| < R $;
		\item $x_1$ lies in the half plane with boundary the    line through $x$ and $x_2$ not containing $b$;
		\item there exists $B^1:=B_R(c_1)$ with $ \{x_1,x\} \subset \partial B^1$, with $\arc(x,x_1)\subset (B^2)^c$, 
		such that the line through $x$ and $x_1$ separates $c_1$ and $x_2$.
	\end{enumerate}
	
	Then the measure of the curved angle $\check{ang}(bxc_1)$ is less than $\frac{\pi}{2}$.
\end{proposition}

\begin{theorem}\label{tangentbound} Let $g$ be an injective $R$-path. 
	Assume that for every $x\in g([a,b])$, $g_x\subseteq  B_R(x)$. 
	
	Then, the measure of $\nor_{\co_R(g_x)}(x)\cap \sfe$ is equal or greater than $\pi/2$.
\end{theorem}

\noindent
{\it Proof.} The proof follows the same lines of the proof of \cite[Theorem 5.2]{MLV2}. The arc length will be replaced by $t$ and $\gamma_x$ by $g_x$; $U^-_x$, $U^+_x$ will be replaced by $\Theta_{-}(x)\cap \sfe$, $\Theta_{+}(x)\cap \sfe$ respectively.

That proof in \cite{MLV2} was divided in three cases: $(a_1)$, $(a_2)$ and $(b)$;
the case $(a_1)$ and the case $(a_2)$ have the same proof except that $(iv)$  of \cite[Proposition 3.2]{MLV2} is replaced by Lemma \ref{Successione}, where $$\overline{u}=\lim_{k \rightarrow +\infty} \frac{g(t^{(k)})-x}{|g(t^{(k)})-x|}.$$
In the case $(b)$ the non intersection property is satisfied by assumption.
\qed

\begin{theorem}\label{corboundtan} Let $g$ be an injective $R$-path. Assume that for every $x\in g([a,b])$, $g_x \subset D(x, \frac{R}{N})$, $N>1$. Then for every $x\in g([a,b])$,
	\begin{equation}\label{PN}
	\meas (\Tan(\co(g_x),x)\cap \sfe) \leq \frac{\pi}{2}+ 2\arcsin \bigg(\frac{1}{2N}\bigg).
	\end{equation}
\end{theorem}

\noindent
{\it Proof.} The proof is the same as in \cite[Theorem 5.4]{MLV2}, with Theorem \ref{tangentbound} instead of \cite[Lemma 4.5]{MLV2}.
\qed

\begin{remark}\label{Wedge gx}
	\textnormal{Let $g$ be an injective $R$-path. Then there is a wedge of opening $\frac{\pi}{2} + 2\arcsin (\frac{1}{2N})< \pi$ containing $g_x$.}
\end{remark}

\begin{definition}\label{Detour}
	Let $f\colon [\alpha,\beta]\to \R^2$ be a rectifiable curve. Let $\alpha\leq\alpha_1< \beta_1\leq\beta$, $y_1=f(\alpha_1)$, $y_2=f(\alpha_2)$ and $y_1 \neq y_2$. The detour of $f$ is defined as
	$$\frac{\length(f_{y_1,y_2})}{|y_1-y_2|}.$$
\end{definition}

\begin{definition}\label{defphiself}
	Let $f\colon [\alpha,\beta]\to \R^2$ be a continuous 
	mapping with the property that for every $y=f(\tau)\in f([\alpha,\beta])$,
	there exists a wedge of opening $\phi\in [0,\pi)$ at the point $y$, which contains $f([\tau, \beta])$; 
	$f$ has been called a $\phi$-self-approaching curve and studied in \cite{Langetepe}.
\end{definition}

\begin{proposition}\label{Langetepe}\cite[Theorem 6]{Langetepe}
	Let $\phi\in [0,\pi)$ and $f$ a $\phi$-self-approaching curve. Then $f$ is a rectifiable curve and
	\begin{equation*}
	(1+\cos(\phi))\length(f)\leq \per(\co(f)).
	\end{equation*}
	Moreover there exists a constant $c({\phi})$, depending on the angle $\phi$ only, for which the detour of $f$ satisfies
	\begin{equation*}
	\frac{\length(f_{y_1,y_2})}{|y_1-y_2|}\leq c({\phi})
	\end{equation*}
	with $y_i=f(\tau_i)$, $\alpha\leq \tau_1<\tau_2\leq\beta$ and $i=1,2$.
\end{proposition}

\begin{corollary}\label{reflection corollary}
	Let $g$ be an injective $R$-path such that $g_x \subseteq D(x,\frac{R}{N})$, $N>1$, for all $x\in g([a,b])$. Let $h(t)=g(-t)$, $t\in [-b,-a]$; then $h$ is a $\phi$-self-approaching curve with $\phi=\frac{\pi}{2}+2\arcsin(\frac{1}{2N})< \pi$.
\end{corollary}

\noindent
{\it Proof.} It is a consequence of Remark \ref{Wedge gx}.
\qed

\begin{theorem}\label{First Main Thm}
	Let $N>1$, $\phi=\frac{\pi}{2}+2\arcsin(\frac{1}{2N})$ and $x_0\in \R^2$. Let $g\colon [a,b]\to \R^2$ be an injective $R$-path, with $g([a,b])\subseteq D(x_0,\frac{R}{2N})$. Then $g$ is a rectifiable curve and
	\begin{equation}\label{Lenght and per}
	\length(g)\leq\ \frac{\per(\co(g))}{1+\cos(\phi)}.
	\end{equation}
	Moreover the detour of $g$ satisfies
	\begin{equation}\label{detour 2}
	\frac{\length(g_{x_1,x_2})}{|x_1-x_2|}\leq c(\phi).
	\end{equation}
\end{theorem}

\noindent
{\it Proof.} Let $h=h(t)$ be the path defined in Corollary \ref{reflection corollary}; $h$ is an injective $\phi$-self-approaching curve. Obviously the detour of $h$ is equal to the detour of $g$. Then by Proposition \ref{Langetepe}, the inequalities \eqref{Lenght and per} and \eqref{detour 2} hold.
\qed

\section{\textbf{All $R$-paths are injectives}}

\begin{definition}
	Let $g\colon [a,b] \to \R^2$ be a continuous function with $g(a)\neq g(b)$. A multiple point $x$ of $g$ is a point such that $\stackrel{-1}{g}(\{ x \})$ is not a single point of $[a,b]$.
\end{definition}

\begin{lemma}\label{Jordancurve}
	Let $x_0\in \R^2$ and $N>1$. Let $g\colon [a,b]\to \R^2$ be an $R$-path with $g(a)\neq g(b)$ and $g([a,b])\subset D(x_0,\frac{R}{2N})$. Then there exists an injective curve $g^{(0)}$ joining $g(a)$ to $g(b)$ that is an $R$-curve and $$g^{(0)}([a,b])\subset g([a,b]).$$
\end{lemma}

\noindent
{\it Proof.}(Partially adapted from \cite{Falconer}) For each multiple point $x$ of $g([a,b])$, let $I_x$ be the largest interval $[t_1,t_2]$ with $g(t_1)=g(t_2)=x$.
There exists a finite or countable collection of closed intervals $I^{(k)}:=I_{x_k}$, $k=1,2,\cdots$, with disjoint interiors such that $g$ is injective outside $\bigcup_{k}\interno(I^{(k)})$. Let us denote $J=\big\{I^{(k)}, \; k = 1,2 \dots \big\} $; a construction of $J$ can be found in the Appendix.

Let $f\colon [a,b]\to [a,b]$ be a continuous, non decreasing, surjective function, with $f(a)=a$, $f(b)=b$, so that if $t_1<t_2$, then $f(t_1)\leq f(t_2)$, with equality if and only if $t_1$ and $t_2$ lie in a some interval of $J$. 


Define $g^{(0)}\colon [a,b]\to \R^2$ by
\begin{equation*}
g^{(0)}(u)=\begin{cases}g(\stackrel{-1}{f}(\{ u \})), & \mbox{if}\ \stackrel{-1}{f}(\{ u \})\ \mbox{is a point of}\ [a,b],\\
x, & \mbox{if}\ \stackrel{-1}{f}(\{ u \})=I_x\in J.
\end{cases}
\end{equation*}

\noindent
Clearly if $a\leq u_1<u_2\leq b$, then $g^{(0)}(u_1)\neq g^{(0)}(u_2)$, so $g^{(0)}$ is an injective curve. Moreover for $u\in [a,b]$ $$g^{(0)}([a,u])\subset g([a,t]),$$ where 

\begin{equation*}
t=\begin{cases}\stackrel{-1}{f}(\{u\}), & \mbox{if}\ \stackrel{-1}{f}(\{u\})\ \mbox{is a point},\\
\min \stackrel{-1}{f}(\{u\}), & \mbox{if}\ \stackrel{-1}{f}(\{u\})\in J.
\end{cases}
\end{equation*}
\noindent
Similarly $g^{(0)}([u,b])\subset g([t,b])$, where 

\begin{equation*}
t=\begin{cases}\stackrel{-1}{f}(\{u\}), & \mbox{if}\ \stackrel{-1}{f}(\{u\})\ \mbox{is a point},\\
\max \stackrel{-1}{f}(\{u\}), & \mbox{if}\ \stackrel{-1}{f}(\{u\})\in J.
\end{cases}
\end{equation*}

\vspace{0.25cm}
\noindent
Let us show that $g^{(0)}$ is an $R$-path. 

Let us consider
\begin{equation*}
\Phi_-=-\Tan(g^{(0)}([a,u]),g^{(0)}(u))\quad \textnormal{and}\quad \Phi_+=\Tan(g^{(0)}([u,b]),g^{(0)}(u)).
\end{equation*}
If $\stackrel{-1}{f}(\{u\})$ is a point $t$, then $g^{(0)}(u)=g(t)$ and
\begin{equation*}
\Phi_{-}\subset\Theta_{-}(g(t))
\end{equation*}
and
\begin{equation*}
\Phi_{+}\subset \Theta_{+}(g(t)).
\end{equation*}
If $l \in \Phi_{-}\cup \Phi_{+}$, with $|l|=1$, then 
$$g^{(0)}([a,u])\subset g([a,t])\subset B^c(g(t)+Rl)=B^c(g^{(0)}(u)+Rl).$$
If $\stackrel{-1}{f}(\{u\})$ is an interval, let $\stackrel{-1}{f}(\{u\})=[\alpha,\beta]\subset [a,b]$. Then 
\begin{equation*}
\Phi_{-}\subset-\Tan(g([a,\alpha]),g(\alpha))
\end{equation*}
and
\begin{equation*}
\Phi_{+}\subset\Tan(g([\beta,b]),g(\beta)).
\end{equation*}
Recall that $g^{(0)}([a,u])\subset g([a,\alpha])$. Let $a\leq v<u$, then $g^{(0)}(v)=g(t)$ for some $t\in [a,\alpha]$ and
\begin{equation*}
B^c(g^{(0)}(u)+Rl)=B^c(g(\alpha)+Rl)\ni g(t)\quad \textnormal{for}\ l\in \Phi_{-}
\end{equation*}
and
\begin{equation*}
B^c(g^{(0)}(u)+Rl)=B^c(g(\beta)+Rl)\ni g(t)\quad \textnormal{for}\ l\in \Phi_{+}.
\end{equation*}
\noindent
In both cases $g^{(0)}(v)\in B^c(g^{(0)}(u)+Rl)$, for $a\leq v\leq u\leq b$, thus $g^{(0)}$ is an $R$-path.

As $g^{(0)}$ is an injective $R$-path, by Theorem \ref{First Main Thm}, $g^{(0)}$ is a rectifiable curve and by Remark \ref{equivalent R-path} is an $R$-curve.
\qed

\begin{theorem}\label{Second Main Thm}
	Let $g\colon [a,b]\to \R^2$, $g([a,b])\subset D(x_0,\frac{R}{2N})$, $N>1$, be an $R$-path.
	Then $g$ can be reparametrized in an injective $R$-path.
\end{theorem}

\noindent
{\it Proof.} Let us reparametrize the path in such way that $\stackrel{-1}{g}(\{x\})$ has no interior points for every $x$, see Remark \ref{Joinying curve}. Let us show now that $g$ has not multiple points.

By contradiction, let us assume that there exist $a\leq \alpha<\beta\leq b$ satisfying $g(\alpha)=g(\beta)$.

\noindent
Since $\stackrel{-1}{g}(\{g(\alpha)\})$ is a closed set without interior points, then there exists a strictly decreasing sequence $\{ \alpha_n \}_{n\in \N}\subset [\alpha,\beta]$ such that
\begin{equation*}
\lim_{n \rightarrow +\infty}\alpha_n=\alpha,\quad  g(\alpha_{n-1})\neq g(\alpha_{n}),\ g(\alpha_n)\neq g(\beta),\ \forall\ n\in \N.
\end{equation*}
{\it Claim.} Let us show that there exists a sequence of $R$-curves $$g_n\colon [\alpha_n,\beta]\to \R^2,\quad n\in \N,$$ satisfying
\begin{equation*}
g_{n-1}([\alpha_{n-1},\beta])\subset g_{n}([\alpha_n,\beta])\subset g([\alpha,\beta]),\quad g_n(\alpha_n)=g(\alpha_n),\ g_n(\beta)=g(\beta),
\end{equation*}
and $g_n([\alpha_n,\beta])\subset D(x_0,\frac{R}{2N})$. 

\vspace{0.25cm}
\noindent
{\it Proof of claim.} By Lemma \ref{Jordancurve}, there exists an injective $R$-path that is also an $R$-curve
$$g^{(1)}\colon [\alpha_1,\beta]\to \R^2$$
satisfying
$$g^{(1)}(\alpha_1)=g(\alpha_1),\quad g^{(1)}(\beta)=g(\beta)\quad \textnormal{and}\quad g^{(1)}([\alpha_1,\beta])\subset g([\alpha_1,\beta]).$$
Assume $g_1=g^{(1)}$.

\noindent
Let $n>1$. Let us assume that $g_{n-1}$ is already defined on $[\alpha_{n-1},\beta]$ and $g_{n-1}$ is an injective $R$-path; let us construct $g_n$. Since $g(\alpha_{n-1})\neq g(\alpha_{n})$ and $g|_{[\alpha_{n},\alpha_{n-1}]}$ is an $R$-path, by Lemma \ref{Jordancurve}, there exists an injective $R$-curve $$g^{(n)}\colon [\alpha_n,\alpha_{n-1}]\to \R^2$$ satisfying
$$g^{(n)}(\alpha_n)=g(\alpha_n),\quad g^{(n)}(\alpha_{n-1})=g(\alpha_{n-1})=g_{n-1}(\alpha_{n-1})$$ and $$g^{(n)}([\alpha_n,\alpha_{n-1}])\subset g([\alpha_n,\alpha_{n-1}]).$$

Now define $g_n$ as the $R$-curve joining $g^{(n)}$ and $g_{n-1}$ according to Lemma \ref{joining curve}.

\noindent
Then $g_n$ is injective, $g_{n-1}([\alpha_{n-1},\beta])\subset g_n([\alpha_n,\beta])\subset g([\alpha_n,\beta])$ and $g_n(\alpha_n)=g(\alpha_n)$, $g_n(\beta)=g(\beta)$. 

The claim is proved.

\vspace{0.25cm}
\noindent
By \eqref{detour 2} of Theorem \ref{First Main Thm}, 
$$\frac{\length(g_n([\alpha_n,\beta]))}{|g_n(\alpha_n)-g_n(\beta)|}\leq c(\phi).$$

\vspace{0.25 cm}
On the other hand $\length(g_n([\alpha_n,\beta]))\geq \length(g_1([\alpha_1,\beta]))\geq |g_1(\alpha_1)-g(\beta)|>0$, thus
$$|g(\alpha_n)-g(\beta)|\geq \frac{\length(g_1([\alpha_1,\beta]))}{c(\phi)}\geq \frac{|g_1(\alpha_1)-g(\beta)|}{c(\phi)}>0 .$$
This is impossible as 
$$\lim_{n \rightarrow +\infty}g(\alpha_n)=g(\alpha)=g(\beta).$$
Contradiction.
\qed

\begin{theorem}\label{R-path and curve}
	Let $g\colon [a,b]\to \R^2$ be an $R$-path. Then $g$ is an $R$-curve.
\end{theorem}

\noindent	
{\it Proof.} Let $N>1$. The compact set $g([a,b])$ is covered by a finite number of open circles $C_i$, $i=1,\cdots,n$, centered in $g([a,b])$ with radius $\frac{R}{2N}$, such that there exist $a=x_0<x_1< \cdots < x_i< \cdots < x_n=b$, with $g([x_i,x_{i+1}])\subset C_i$.

The map $g$ restricted to $[x_i,x_{i+1}]$ is injective and rectifiable by Theorem \ref{Second Main Thm} and by Theorem \ref{First Main Thm}.
Then $g$ is a rectifiable curve and, by Remark \ref{equivalent R-path}, it is an $R$-curve.
\qed

\begin{remark}($R$-curves and $\lambda$-eels)
	\textnormal{Let us show that the families of $R$-curves and $\lambda$-eels are different.}
	
	\textnormal{In \cite[Example 1]{steepest I} it is shown that there exists an $R$-curve that is not a $\phi$-self-approaching curve (defined in Definition \ref{defphiself}). It has been noticed in \cite[Remark 6.1 and Lemma 6.2]{Daniilidis/Deville} that a $\lambda$-eel $\gamma$ is a $\phi$-self-approaching curve for some $\phi$. Thus there are $R$-curves that are not $\lambda$-eels.}
	
	\textnormal{On the other hand in \cite[Example 2]{steepest I} it is introduced a $\lambda$-eel (for suitable $\lambda$) that is not an $R$-curve for every $\lambda>0$.}
\end{remark}

\begin{theorem}\label{R-path in a disk}
	Let $g$ be an $R$-path with image contained in an open circle of radius $R$ and representation $x=x(s)$, $0\leq s\leq L$, in arc length.
	
	Then there exists a nested family $\Omega_s$ of sets of reach $\geq R$ with the property that $x(s)\in \partial \Omega_s$ and $x'(s)$ belongs to the normal cone to $\Omega_s$ at $x(s)$.
\end{theorem}	

\noindent
{\it Proof.} Let us define $$\Omega_s=\co_R(x([0,s])).$$
By Remark \ref{R-hull inclusion}, if $s_1<s_2$, then $$\co_R(x([0,s_1]))\subset \co_R(x([0,s_2])).$$
Thus $\Omega_s$ is a nested family of sets of reach $\geq R$, such that $x(s)\in \partial \Omega_s$ and $x'(s)\in \nor(\Omega_s,x(s))$.
\qed

\vspace{0.50cm}
\noindent\textbf{Acknowledgements.}
	This work has been partially supported by INDAM-GNAMPA(2019).
	
	This work is dedicated to our friend Sergio Vessella.

\section*{\textbf{Appendix}}

Here a proof is given of the following fact quoted in \cite{Falconer}.

\begin{theorem*}
	Let $g\colon [a,b]\to \R^2$ be a continuous mapping, $g(b)\neq g(a)$. There exists a finite or countable collection $J=\big\{I^{(k)}, \; k = 1,2 \dots \big\} $ of proper closed intervals in $[a,b]$, with disjoint interiors, so that $g$ is injective outside of $\bigcup_{k}\interno(I^{(k)})$ and, if $[\alpha,\beta]\in J$, then $g(\alpha)=g(\beta)$.
\end{theorem*}

\noindent
{\it Proof.} Let $K\subset [a,b]$ be a compact subset. Let us define
\begin{equation*}
\phi_K(t):=\max \{ \tau\in K|\ g(\tau)=g(t) \}.
\end{equation*}

\noindent
{\it Claim.} The following properties hold:
\renewcommand{\theenumi}{\roman{enumi}}
\begin{enumerate}
	\item $\forall\ t\in K$, $t\leq \phi_K(t)$ and  $\phi_K(t)-t\leq b-a$.
	\item If $\phi_K(t)-t\equiv 0$ for a compact set $K\subset [a,b]$, then $g_{|_K}$ is injective.
	\item $\phi_K(t)-t$ has a maximum in $K$.
\end{enumerate}

\noindent
{\it Proof of claim.} The first two statements are obvious, let us prove \textit{(iii)}.

Let $L=\sup_{t\in K }\{ \phi_K(t)-t \}\in [0,b-a]$. If $K$ is finite, then the assertion is obvious.

\noindent
If $K$ is not finite, then there exists a sequence $\{ t_n\}_{n\in \N}\subseteq K$ such that $\phi_K(t_n)-t_n$ is increasing and $$\lim_{n \rightarrow +\infty} \phi_K(t_n)-t_n=L.$$
\noindent
Let $\{ t_{n_j}\}_{j\in \N}$ be a subsequence of $\{ t_n\}_{n\in \N}$ satisfying
\begin{equation*}
\lim_{j \rightarrow +\infty}t_{n_j}=\overline{t}\in K,\quad \lim_{j \rightarrow +\infty}\phi_K(t_{n_j})=\overline{\phi}\in K.
\end{equation*}
As $g(t_{n_j})=g(\phi_K(t_{n_j}))$, then $g(\overline{t})=g(\overline{\phi})$. Thus $\overline{\phi} \leq \phi_K(\overline{t})$ and $$L=\overline{\phi}-\overline{t}\leq \phi_K(\overline{t})-\overline{t}\leq L,$$
thus $\overline{\phi}=\phi_K(\overline{t})$ and $L=\phi_K(\overline{t})-\overline{t}$ is a maximum.

The claim is proved.

\vspace{0.50cm}

For each compact subset $K$ of $[a,b]$, let $M(K)$ be the maximum in $K$ of $\phi_K(t)-t$, $\mu(K)$ be the Lebesgue measure of $K$ and $\delta(K)$ be the maximum of the length of the closed subintervals of $K$; $\delta(K)=0$ if in $K$ there are no intervals.

\vspace{0.15cm}
\noindent
{\it Construction of $J$.}
\noindent

In what follows, when $ I $ is a closed interval, $ \interno(I)  $ will be denoted by $\mathcal{I}$.
Let us construct a sequence of nested compact sets 
\begin{equation*}
K_0=[a,b]\supset K_1 \supset K_2 \supset \cdots 
\end{equation*}
with $K_n$ given by $[a,b]$ minus a finite number of open disjoint intervals, where $g$ has the same value at their extreme points of each intervals:
$$ K_{n} = [a,b] \setminus \cup_{l=1}^{s_n}  \mathcal{I}_{l} ^{(n)}.
$$

\noindent
Let us assume that $K_{n-1}$ has been constructed, with $\delta(K_{n-1})>0$ and $\mu(K_{n-1})<b-a$; let us define $K_n$ and let us show that $\delta(K_n)>0$.

\vspace{0.15cm}
\noindent
In $K_{n-1}$ there is a family of open intervals of measure $M(K_{n-1})$; these intervals may overlap; let $\mathcal{I}_{1}^{(n-1)}$ be the leftmost one.

\vspace{0.15cm}
\noindent
Let us define $K_n^{(1)}:=K_{n-1}\setminus \mathcal{I}_{1}^{(n-1)}$; then $K_n^{(1)}$ is $[a,b]$ minus a finite number of open disjoint intervals and 
$$\delta_n^{(1)}:=\delta(K_n^{(1)})>0,\quad \mu_n^{(1)}:=\mu_n(K_n^{(1)})=\mu(K_{n-1})-M(K_{n-1})$$ 
hold.

If $M(K_n^{(1)})<M(K_{n-1})$, then $K_n:=K_n^{(1)}$, $\mu(K_n)=\mu(K_n^{(1)})$, $0<\delta(K_n)=\delta_n^{(1)}\leq \delta(K_{n-1})$.

\vspace{0.25cm}
\noindent
If	$M(K_n^{(1)})=M(K_{n-1})$, then let us consider the open disjoint intervals in $K_n^{(1)}$ with measure $M(K_{n-1})$; these intervals may overlap; let $\mathcal{I}_{2}^{(n-1)}$ be the leftmost one.

\vspace{0.15cm}
\noindent
Let us define $K_n^{(2)}:=K_n^{(1)}\setminus \mathcal{I}_{2}^{(n-1)}$; then $K_n^{(2)}$ is $[a,b]$ minus a finite number of open disjoint intervals and 
$$\delta_n^{(2)}:=\delta(K_n^{(2)})>0,\quad \mu_n^{(2)}:=\mu_n^{(1)}-M(K_{n-1})=\mu(K_{n-1})-2M(K_{n-1})>0$$ hold.

If $M(K_n^{(2)})<M(K_{n-1})$, then $K_n:=K_n^{(2)}$; otherwise the procedure have to be iterated by constructing $K_n^{(j)}$, $j\geq 3$. Let us notice that $j$ satisfies  
$$\mu_n^{(j)}=\mu(K_{n-1})-jM(K_{n-1})>0,\quad \forall\ j.$$

\vspace{0.15cm}
\noindent
Thus $K_n$ is constructed in a finite number of iterations and $\delta(K_n)>0$. The compact set $K_n$ is given by $[a,b]$ minus a finite number of open disjoint intervals, where $g$ has the same value at their extreme points; $K_n$ is the union of a finite non zero number of closed intervals (where $g$ does not have the same values at their extreme points) and a finite number of isolated points; moreover, if $ \mathcal{I} $ is an interval in $ [a,b]\setminus K_{n-1}$, then $ \mathcal{I} \subset  [a,b]\setminus K_{n}. $

\vspace{0.25cm}
\noindent
Let $M_n=M(K_n)$. There are two possibilities.
\renewcommand{\theenumi}{\roman{enumi}}
\begin{enumerate}
	\item If $M_{n_{0}}=0$ for some $n_{0}\in \N$, then $g_{|_{K_{n_{0}}}}$ is injective where 
	\begin{equation*}
	K_{n_{0}}=[a,b]\setminus \cup_{l=1}^{s_{n_{0}}}  \mathcal{I}_{l} ^{(n_{0})} 
	\end{equation*}
	and $J$ can be defined as 
	
	$$
	J=\big\{ cl \; \mathcal{I}^{n_{0}}_{l}: l=1, \dots s_{n_{0}}  \big\}.$$
	
	\item If $M_n > 0$, for every $n$, then the sequence $M_0,\cdots,M_n$ is decreasing. Moreover we have
	$$M_0+\cdots +M_{n-1}\leq b-a,\quad \lim_{n \rightarrow +\infty} M_n=0$$
	and  $J$ can be defined as  	
	
	$$
	J=\big\{ cl \; \mathcal{I}^{n}_{l}: l=1, \dots s_{n}, n \in \N  \big\},$$
	where $\mathcal{I}^{n}_{l}$ appears infinite times.

\end{enumerate}

In both cases  $J$ is a finite or countable collection of closed intervals and it satisfies all the stated properties.
\qed

\bigskip

N. Lombardi:\\ 
Vienna University of Technology, Institute of Discrete Mathematics and Geometry,\\
Wiedner Hauptstra\ss e 8--10/104, 1040 Wien, Austria.

\vspace{0.25cm}

M. Longinetti:\\
DIMAI, Universit\`a  di Firenze, V.le Morgagni 67/a, 50134, Italy.\\
Electronic mail addresses: marco.longinetti@unifi.it
 
\vspace{0.25cm}
P. Manselli:\\
DIMAI, Universit\`a  di Firenze, V.le Morgagni 67/a, 50134, Italy.

\vspace{0.25cm}
A. Venturi:\\
Universit\`a  di Firenze, P.le delle Cascine 15, 50144, Italy.


\begin{thebibliography}{}
	
	\bibitem{Daniilidis/Deville} A. Daniilidis, R. Deville, E. Durand-Cartagena, {Metric and geometric ralaxations of self-contracted curves}, Journal of Optimization Theory and Applications 182, 81-109 (2019).
	
	\bibitem{Mainik} I.F. Ma\v{i}nik, {An estimate of the length of the curves of descent}, Sibirsk. Mat. Zh. {\bf 33}, 215-218 (1992).
	
	\bibitem{Manselli-Pucci} P. Manselli, C. Pucci,
	{Maximum length of Steepest descent curves for quasi-convex Functions}, Geometriae Dedicata {\bf 38}, 211-227 (1991).
	
	\bibitem{MLV} M. Longinetti, P. Manselli, A. Venturi, {On  steepest descent curves for quasi convex families in $\R^n$}, Math. Nachr. {\bf 288}, 420-442 (2015). 
	
	\bibitem{Langetepe} C. Icking, R. Klein and E. Langetepe, {Self-approaching curves}, Math. Proceedings Cambridge Philos. Sc. {\bf 125}, 441-453 (1999).
	
	\bibitem{Daniilidis} A. Daniilidis, O. Ley, S. Sabourau, {Asymptotic behavior of self-contracting planar curves and gradient orbits of convex functions}. J Math Pures Appl. 2010;94:183-199.
	
	\bibitem{Langetepe2} 0. Aichholzer, F. Aurenhammer, C. Icking, R. Klein, E. Langetepe and  G. Rote,
	{Generalized self-approaching curves}, Discr. Appl. Math.{\bf 109}, 3-24 (2001).
	
	\bibitem{Mainik2} I.F. Ma\v{i}nik, {On stability in a theorem on an estimate for the length of a descent curve}, Sibirsk. Mat. Zh. {\bf 38}, 867-875 (1997).
	
	\bibitem{MLV2} M. Longinetti, P. Manselli, A. Venturi, {Plane $R$-curves II}, https://doi.org/10.1080/00036811.2019.1600677.
	
	\bibitem{steepest I} M. Longinetti, P. Manselli, A. Venturi, {Plane $R$-curves and their steepest descent properties I}, Appl. Anal. {\bf 98}, (2018).
	
	\bibitem{Federer} H. Federer, { Curvature measures}, Trans. Amer. Math. Soc. {\bf 93}, 418-481(1959).
	
	\bibitem{Ratay1} J. Rataj, {Determination of spherical area measures by means of dilatation volumes}, Math. Nachr. {\bf 235}, 143-162 (2002). 
	
	\bibitem{Colman} A. Colesanti, P. Manselli, {Geometric and Isoperimetric Properties of sets of Positive Reach in $\E^d$ }, Atti Semin. Mat. Fis. Univ. Modena Reggio Emilia {\bf 57}, 97-113 (2010).
	
	\bibitem{Falconer} K.T. Falconer, The geometry of fractal sets, Cambridge University Press, 1985.
	
	
	
	
	
	
	
	
	
	
	
	
	
	
	
	
\end{thebibliography}
\end{document}